\documentclass{article}     
\usepackage[utf8]{inputenc}
\usepackage{amsfonts}
\usepackage[english]{babel}
\usepackage{xcolor}
\usepackage[T1]{fontenc}
\usepackage{amsmath, amssymb,epic,graphicx,mathrsfs,enumerate}
\usepackage[all]{xy}
\usepackage{color}
\usepackage{comment}
\usepackage{enumitem}
\usepackage{esint}
\usepackage{hyperref}
\usepackage{amsthm}
\usepackage{latexsym}
\usepackage{epsfig}
\usepackage{soul}
\usepackage{geometry}
\usepackage{doi}
\hypersetup{colorlinks=true,
            linkcolor=black,
            filecolor=black,
            urlcolor=blue,
            citecolor=blue}

\newcommand{\R}{\mathbb{R}}
\newcommand{\e}{\varepsilon}

\newcommand{\mail}[1]{\href{mailto:#1}{\texttt{#1}}}

\newtheorem{theorem}{Theorem}[section]
\newtheorem{corollary}[theorem]{Corollary}
\newtheorem{lemma}[theorem]{Lemma}
\newtheorem{prop}[theorem]{Proposition}
\newtheorem{definition}{Definition}[section]
\theoremstyle{remark}
\newtheorem{obs}[theorem]{Remark}

\title{Existence and asymptotic behavior of non-normal conformal metrics on $\mathbb{R}^4$ with sign-changing $Q$-curvature}

\author{Chiara Bernardini}

\begin{document}

\maketitle

\begin{abstract}
\noindent We consider the following prescribed $Q$-curvature problem 
\begin{equation}\label{uno}
\begin{cases}
\Delta^2 u=(1-|x|^p)e^{4u}, \quad\text{on}\,\,\mathbb{R}^4\\
\Lambda:=\int_{\mathbb{R}^4}(1-|x|^p)e^{4u}dx<\infty.
\end{cases}
\end{equation}
We show that for every polynomial $P$ of degree 2 such that $\lim\limits_{|x|\to+\infty}P=-\infty$, and for every $\Lambda\in(0,\Lambda_\mathrm{sph})$, there exists at least one solution to problem \eqref{uno} which assume the form $u=w+P$, where $w$ behaves logarithmically at infinity. Conversely, we prove that all solutions to \eqref{uno} have the form $v+P$, where 
$$v(x)=\frac{1}{8\pi^2}\int\limits_{\mathbb{R}^4}\log\left(\frac{|y|}{|x-y|}\right)(1-|y|^p)e^{4u}dy$$
and $P$ is a polynomial of degree at most 2 bounded from above. Moreover, if $u$ is a solution to \eqref{uno}, it has the following asymptotic behavior
$$u(x)=-\frac{\Lambda}{8\pi^2}\log|x|+P+o(\log|x|),\quad\text{as}\,\,|x|\to+\infty.$$
As a consequence, we give a geometric characterization of solutions in terms of the scalar curvature at infinity of the associated conformal metric $e^{2u}|dx|^2$.\\



\end{abstract}

\section{Introduction}

Let $M$ be a 4-dimensional Riemannian manifold endowed with a metric $g$. If $\mathrm{Ric}_g$ denotes the Ricci tensor of $(M,g)$, $R_g$ the scalar curvature and $\Delta_g$ the Laplace-Beltrami operator, the $Q$-curvature $Q_g^4$ and the Paneitz operator $P_g^4$ are defined as follows
$$Q_g^4:=-\frac{1}{6}\left(\Delta_g R_g-R_g^2+3|\mathrm{Ric}_g|^2\right)$$
$$P_g^4(f):=\Delta_g^2 f+\mathrm{div}\left(\frac{2}{3}R_g g-2\mathrm{Ric}_g \right)df,\quad\text{for}\,\,f\in C^\infty(M).$$
These definitions have been introduced by T. Branson and B. {\O}rsted \cite{BO} and S. Paneitz \cite{P}, and then generalized to higher order $Q$-curvatures $Q_g^{2m}$ and Paneitz operators $P_g^{2m}$ on a $2m$-dimensional manifold (see e.g. \cite{B,B1,FG,FH,GJMS}). The Paneitz operator is a sort of higher order Laplace-Beltrami operator and the $Q$-curvature can be seen as the higher order counterpart of the Gaussian curvature, this is also pointed out by the fact that in dimension 2 we have $P_g^2=-\Delta_g$ and $Q_g^2=K_g$.
The study of Paneitz operator and Q-curvature has gained a lot of attention in conformal geometry due to their covariant properties. The total $Q$-curvature is a global conformal invariant, namely if $M$ is closed and $g_u=e^{2u}g$, we have $$\int_M Q^{2m}_{g_u}d\,\mathrm{vol}_{g_u}=\int_M Q^{2m}_g d\, \mathrm{vol}_{g}$$
moreover, this integral gives information on the topology of the manifold (see the Gauss-Bonnet-Chen's Theorem \cite{Ch}). 
Recall also that the Paneitz operator is conformally invariant, this means that if we consider the conformal metric $g_u=e^{2u}g$, we have 
$$P_{g_u}^{2m}=e^{-2mu}P_g^{2m}$$
and it satisfies a generalised version of the Gauss equation
$$P_g^{2m}u+Q_g^{2m}=Q_{g_u}^{2m}e^{2mu}.$$
On $\mathbb{R}^{2m}$ with the standard Euclidean metric $|dx|^2$ we have 
$$P^{2m}_{|dx|^2}=(-\Delta)^m \quad \text{and} \quad Q^{2m}_{|dx|^2}\equiv0.$$

A classical problem in conformal geometry is the ``problem of prescribing $Q$-curvature'', that is finding whether a given function $K$ on a manifold $(M,g)$, can be the $Q$-curvature of a conformal metric $g_u=e^{2u}g$. Therefore, the question concerns the set of solutions to the following equation 
\begin{equation}\label{1a}
P^{2m}_{g} u+Q^{2m}_{g}=Ke^{2mu}
\end{equation}
where $K=Q^{2m}_{g_u}$ is the \textbf{prescribed} $Q$-curvature of the metric  $g_u$. If we consider the Euclidean space $\mathbb{R}^{2m}$ endowed with the standard Euclidean metric $|dx|^2$, equation \eqref{1a} becomes
\begin{equation}\label{1b}
(-\Delta)^m u = K e^{2mu}, \quad\text{in}\,\,\mathbb{R}^{2m}.
\end{equation}
Due to its geometric meaning, the constant $Q$-curvature case for equation \eqref{1b} has been extensively studied (see e.g. \cite{CC,CL1,CY,HM14, Lin,Mclass,WY} and the references within). The problem which we address, namely the case of non-normal metrics with sign-changing prescribed $Q$-curvature having power-like growth, is (to the best of our knowledge) still unexplored.\\


Let $p>0$ be fixed, we consider the  prescribed $Q$-curvature equation 
\begin{equation}\label{1}
\Delta^2 u=(1-|x|^p)e^{4u}, \quad\text{in}\,\,\mathbb{R}^4
\end{equation}
under the assumption
\begin{equation}\label{2}
\Lambda:=\int_{\mathbb{R}^4}(1-|x|^p)e^{4u}dx<\infty.
\end{equation}
Geometrically, this means that if $u$ is a solution to \eqref{1}-\eqref{2}, then the metric $g_u:=e^{2u}|dx|^2$, which is conformal to the Euclidean metric on $\mathbb{R}^4$, has $Q$-curvature equal to $1-|x|^p$ and finite total $Q$-curvature $\Lambda$. In what follows, we will always assume that $(1-|x|^p)e^{4u}\in L^1(\mathbb{R} ^4)$. Let $u$ be a solution to \eqref{1}-\eqref{2}, we define $v$ as
\begin{equation}\label{vdef}
v(x):=\frac{1}{8\pi^2}\int\limits_{\mathbb{R}^4}\log\left(\frac{|y|}{|x-y|}\right)(1-|y|^p)e^{4u}dy.   
\end{equation}

\begin{definition}[Normal and non-normal solutions]
We call a solution $u$ to \eqref{1}-\eqref{2} normal if there exists a constant $c\in\mathbb{R}$ such that $u$ solves the following integral equation $$u(x)=\frac{1}{8\pi^2}\int\limits_{\mathbb{R}^4}\log\left(\frac{|y|}{|x-y|}\right)(1-|y|^p)e^{4u}dy+c.$$
All other solutions to problem \eqref{1}-\eqref{2} are called non-normal.
\end{definition}

For what concerns normal solutions to \eqref{1}-\eqref{2}, recently, A. Hyder and L. Martinazzi \cite{HM20} proved some existence and non-existence results. In particular, among other things, they showed that problem \eqref{1}-\eqref{2} admits normal solutions if and only if $p\in(0,4)$ and $\Lambda_{*,p}\le\Lambda<\Lambda_{\mathrm{sph}}$ where
$\Lambda_{*,p}:=\left(1+\frac{p}{4}\right)8\pi^2$ and $\Lambda_\mathrm{sph}:=16\pi^2$. Moreover, every normal solution has the following asymptotic behavior 
$$u(x)=-\frac{\Lambda}{8\pi^2}\log|x|+C+O(|x|^{-\alpha}), \quad\text{as}\,\,|x|\to\infty$$
for every $\alpha\in[0,1]$ such that 
$\alpha<\frac{\Lambda-\Lambda_{*,p}}{2\pi^2 }$.

Motivated by the above results, we study the properties of more general solutions (not necessarily normal) to problem \eqref{1}-\eqref{2}. 
Although for $p\ge4$, problem \eqref{1}-\eqref{2} admits no normal solutions, we prove that non-normal solutions do exist. To this end, we consider all polynomials $P$ of degree $2$ such that
$$\lim\limits_{|x|\to\infty}P(x)=-\infty$$
and define the set 
$$\mathcal{P}_2:=\{P\,\,\text{polynomial in}\,\, \mathbb{R}^{4}\,|\,\mathrm{deg}P=2\quad\text{and}\quad \lim\limits_{|x|\to\infty}P(x)=-\infty \}.$$
By means of a result of S.-Y. A. Chang and W. Chen (see Theorem 2.1 in \cite{CC} where, under suitable assumptions on the curvature $K$, using a variational approach, they prove existence of at least one solution to equation $(-\Delta )^{\frac{n}{2}}u=K(x)e^{nu}$ in $\mathbb{R}^{n}$) we prove the following
theorem, which extends to the non-normal case existence results in \cite{HM20}.

\begin{theorem}\label{1.1}
Let $p>0$ be fixed. Then for any $P\in \mathcal{P}_2$ and for every $\Lambda\in(0,\Lambda_ \mathrm{sph})$, there exists at least one solution to problem \eqref{1}-\eqref{2} of the form $$u=w+P$$ where
$$w(x)=-\frac{\Lambda}{8\pi^2}\log|x|+C+o(1), \quad\text{as}\,\,|x|\to+\infty.$$
\end{theorem}

\noindent We shall prove the following classification result.
 
\begin{theorem}\label{1.2}
Let $u$ be a solution to \eqref{1}-\eqref{2} such that $(1-|y|^p)e^{4u}\in L^1(\R^4)$ and $v$ as defined in \eqref{vdef}.
Then there exists an upper-bounded polynomial $P$ of degree at most 2 such that $$u=v+P.$$
Moreover, $u$ has the following asymptotic behavior 
\begin{equation}\label{4}
u(x)=-\frac{\Lambda}{8\pi^2}\log|x|+P(x)+o(\log|x|),\quad\text{as}\,\,|x|\to\infty.
\end{equation}
\end{theorem}

Note that $P$ being upper bounded means that $P$ has even degree, and since $P$ has degree at most 2, this implies that $P$ could only have degree 2 or be constant. For this reason, we can rephrase Theorem \ref{1.2} by saying that all solutions to problem \eqref{1}-\eqref{2} have the form $v+P$, where $v$ behaves logarithmically at infinity and $P$ is an upper bounded polynomial of degree 2 if the solution is non-normal, whereas $P$ is constant if the solution is normal. 

\begin{obs}
We can observe that the function $w$ of Theorem \ref{1.1} and the function $v$ of Theorem \ref{1.2} differ by a constant.
\end{obs}

In order to prove Theorem \ref{1.2}, we obtain some suitable upper and lower estimates for the function $v$ (see Lemma \ref{v<} and Proposition \ref{v>new} below). First, we prove the lower estimate \eqref{v>}. To this end, we
need to overcome some difficulties compared to previous works, due to the fact that the $Q$-curvature is not constant and it changes sign (compare e.g. to the proof of Lemma 2.4 in \cite{Lin}). The lower estimate \eqref{v>} will be crucial to obtain a Liouville-type theorem (see Theorem \ref{teoliouville}), also in this case the proof is quite delicate because the estimate for $-v$ contains the singular integral
$$A(x):=\frac{1}{8\pi^2}\int_{B_1(x)}\log(|x-y|)(1-|y|^p) e^{4u}dy.$$
Then we show that the polynomial $P$ is upper-bounded (see Proposition \ref{sup p}). To do this, we take advantage of a useful result of Gorin and the fact that $|A(x)|\le C$. In order to estimate the singular integral $A(x)$ and prove Proposition \ref{v>new}, we take some ideas from the proof of Lemma 13 in \cite{Mclass} and the one of Lemma 2.4 in \cite{Lin}, but our case is more challenging. This is due to the fact that the singular integral $A(x)$ is over $B_1(x)$ and we need a radius $\tau\in(0,1)$ which can be fixed later. In addition, we lack of a good estimate for $\int_{\mathbb{R}^4 \setminus B_R}v^+\,dx$, which indeed is necessary in \cite{Mclass}, and we do not know a priori the sign of $\Delta u$, which is fundamental to apply an Harnack inequality as in \cite{Lin}.

Finally, we characterize solutions in terms of the behavior at infinity of the scalar curvature of the associated conformal metric. 

\begin{theorem}\label{1.3}
Let $u$ be a solution to problem \eqref{1}-\eqref{2} and set $g_u=e^{2u}|dx|^2$. If $u$ is a normal solution, then 
$$\lim\limits_{|x|\to+\infty}R_{g_u}=0.$$
If $u$ is a non-normal solution, we get
$$\liminf\limits_{|x|\to+\infty}R_{g_u}=-\infty.$$
\end{theorem}

Our result proves that normal solutions have a quite different geometry compared to that of non-normal solutions. Moreover, it also points out that, for this problem, the $Q$-curvature and the scalar curvature are independent of each other.\\

\textbf{Open problem} Can we find non-normal solutions to \eqref{1}-\eqref{2} with arbitrarily large but finite total $Q$-curvature $\Lambda$? In the constant $Q$-curvature case C.-S. Lin \cite{Lin} (see also \cite{CC}) proved that all solutions to $$\begin{cases}\Delta^2 u=6e^{4u} \quad \text{in}\,\,\mathbb{R}^4\\
e^{4u}\in L^1(\mathbb{R}^4)\end{cases}$$
satisfy $V\in(0,\text{vol}(S^4)]$, but his approach is no longer applicable to our case, in fact when $u=v+P$ with $P$ polynomial of degree 2 we need $Q(x)=(1-|x|^p)e^{4P}$ to be radially decreasing, which is not true in general. Since $x\cdot \nabla Q(x)$ does not have a fixed sign, using methods from \cite{M13} or \cite{H}, it would be interesting to see whether there exist solutions to \eqref{1}-\eqref{2} with total $Q$-curvature $\Lambda\ge\Lambda_\text{sph}$.


\section{Existence of solutions}

In this section, we take advantage of a result of A. Chang and W. Chen (see Theorem 2.1 in \cite{CC}). Using a variational approach in a Sobolev space defined on a conical singular manifold, they prove existence of at least one solution to equation 
$$(-\Delta)^{\frac{n}{2}}w=K(x)e^{nw}\quad\text{in}\,\,\mathbb{R}^n,$$
in even dimensions, assuming that $K$ is positive somewhere and for some $s>0$, $K(x)=O\left(\frac{1}{|x|^s}\right)$ near infinity. 

\proof[Proof of Theorem \ref{1.1}]
Let us fix $P\in\mathcal{P}_2$ and $\Lambda\in(0,\Lambda_ \text{sph})$. By Theorem 2.1 in \cite{CC} and its proof (refer also to Section 7 in \cite{HM20}) setting $K(x):=(1-|x|^p)e^{4P}$ and $\mu:=1-\frac{\Lambda}{\Lambda_ \mathrm{sph}}\in(0,1)$, one can find at least one solution $w$ to equation
\begin{equation}\label{b}
\Delta^2 w=K(x)\,e^{4w},\quad\text{in}\,\,\mathbb{R}^4
\end{equation}
such that 
$$\int\limits_{\mathbb{R}^4} K(x)\,e^{4w}\,dx=(1-\mu)\Lambda_{\mathrm{sph}}= \Lambda.$$

It follows immediately that $$u:=w+P$$ is the desired solution to problem \eqref{1}-\eqref{2}. 
More precisely, $w$ is of the form 
$$w=\Tilde{w}\circ\Pi^{-1}+(1-\mu)w_0$$
where $w_0=\log\left(\frac{2}{1+|x|^2}\right)$, $\Pi:S^4\to\mathbb{R}^4$ denotes the stereographic projection, $\tilde{w}=\Bar{w}+C$, where $\Bar{w}$ minimize
a certain functional which takes values in a Sobolev space defined on a conical singular manifold, and $C$ is a suitable constant such that 
$$\int K(x)e^{4\tilde{w}}\,dV=(1-\mu)\Lambda_\mathrm{sph}$$
where the corresponding volume element is $dV=e^{4(1-\mu)w_0}dx$. We have
$$P^4_{g_0}\tilde{w}+6(1-\mu)=(K\circ\Pi)e^{-4\mu(w_0\circ \Pi)}e^{4\tilde{w}},$$
from this identity, with the same argument as in Section 7 of \cite{HM20}, we obtain $\tilde{w}\in C^{3,\alpha}(S^4)$ for $\alpha\in(0,1)$. In particular, $\tilde{w}$ in continuous at the South pole $S=(0,0,0,0,-1)$, which implies
$$w(x)=(1-\mu)w_0(x)+\tilde{w}(S)+o(1)=-\frac{\Lambda}{8\pi^2}\log|x|+C+o(1),\quad\text{as}\,\,|x|\to\infty.$$
\endproof

\begin{obs} 
If $P\in \mathcal{P}_2$ is a radially symmetric polynomial, there exists at least one non-normal radial solution to problem \eqref{1}-\eqref{2} of the form $u=w+P$. This follows from the fact that we can minimize the previous functional over radial functions and obtain $\Bar{w}$ radially symmetric.
\end{obs}


\section{Asymptotic behavior}

In all this section, let $u$ be a solution to problem \eqref{1}-\eqref{2}, we define 
\begin{equation}\label{v}
v(x)=\frac{1}{8\pi^2}\int\limits_{\mathbb{R}^4}\log\left(\frac{|y|}{|x-y|}\right)(1-|y|^p)e^{4u}dy.
\end{equation}
Obviously, we have $\Delta^2 v(x)=(1-|x|^p)e^{4u}$ in $\mathbb{R}^4$.

\begin{lemma}\label{v<}
For $|x|\ge4$, there exists a constant $C$ such that
\begin{equation}\label{valto}
v(x)\le-\frac{\Lambda}{8\pi^2}\log|x|+C.
\end{equation}
\end{lemma}

\proof
For $|x|\ge4$, we decompose $\R^4=B_1(0)\cup A_1\cup A_2\cup A_3$ where
$A_1=\{y\,|\,|y-x|\le|x|/2\}$, $A_2=\{y\,|\,\,1\le|y|\le2\}$ and $A_3=\mathbb{R}^4\setminus(A_1\cup A_2\cup B_1)$. For $y\in B_1$ we have $$\log\left(\frac{|y|}{|x-y|}\right)\le-\log|x|+C$$
hence
$$\int_{B_1}\log\left(\frac{|y|}{|x-y|}\right)(1-|y|^p)e^{4u}dy\le(-\log|x|+C)\int_{B_1}(1-|y|^p)e^{4u}dy.$$
For $y\in A_1$ we have $\log\left(\frac{|y|}{|x-y|}\right)\ge0$ hence the integral over $A_1$ is negative. 
For what concern $A_2$ we have
\begin{align*}
\int_{A_2}\log\left(\frac{|y|}{|x-y|}\right)(1-&|y|^p)e^{4u}dy=\int_{A_2}\log(|y|)\,(1-|y|^p)e^{4u}dy-\int_{A_2}\log(|x-y|)(1-|y|^p)e^{4u}dy\\
\intertext{using the fact that the first integral is non positive, we get}
&\le -\int_{A_2}\log(|x-y|)(1-|y|^p)e^{4u}dy\le -\log|x|\int_{A_2}(1-|y|^p)e^{4u}dy+C
\end{align*}
where in the last inequality we used that for $y\in A_2$ $\log|x-y|\le\log|x|+C$.
For $y\in A_3$ since $|x-y|\le |x|+|y|\le |x|\,|y|$ we have 
$$\log\left(\frac{|y|}{|x-y|}\right)\ge-\log|x|$$
and since in this case $1-|y|^p<0$ we get
$$\log\left(\frac{|y|}{|x-y|}\right)(1-|y|^p)\le-\log(|x|)\,(1-|y|^p)$$
and hence
$$\int_{A_3}\log\left(\frac{|y|}{|x-y|}\right)(1-|y|^p)e^{4u}dy\le-\log(|x|)\int_{A_3}(1-|y|^p)e^{4u}dy.$$
Summing up, we finally obtain
$$v(x)\le-\frac{1}{8\pi^2}\log(|x|)\int_{A_1^C}(1-|y|^p)e^{4u}dy+C,$$
since $\int_{A_1^C}(1-|y|^p)e^{4u}dy\ge\Lambda$ we have
$$v(x)\le-\frac{\Lambda}{8\pi^2}\log|x|+C.$$
\endproof

\begin{lemma}\label{v>=}
For any $\varepsilon>0$, there exists $R=R(\varepsilon)>0$ such that for $|x|\ge R$ 
\begin{equation}\label{v>}
v(x)\ge -\frac{1}{8\pi^2}(\Lambda+5\varepsilon)\log|x|-\frac{1}{8\pi^2}\int_{B_1(x)}\log(|x-y|)(1-|y|^p) e^{4u}dy.
\end{equation}
\end{lemma}

\proof
We can choose $R_0=R_0(\e)>1$ such that
$$\int_{B_{R_0}}(1-|y|^p)e^{4u}dy\le\Lambda+\e.$$
Let us take $R>2R_0$ and assume that $|x|\ge R$, we can decompose
$$\mathbb{R}^4=B_{R_0}(0)\cup A_1 \cup A_2$$
where
$$A_1:=\{y\in\R^4\,:\,|y-x|\le|x|/2\}$$
$$A_2:=\{y\in\R^4\,:\,|y-x|>|x|/2,\,\,|y|\ge R_0\}.$$
For $|x|\ge R$ and $|y|\le R_0$ we have
$$\log\left(\frac{|y|}{|x-y|}\right)\le-\log|x|+C<0$$
hence, we get
\begin{align*}
\int_{B_{R_0}\setminus B_1}\log\left(\frac{|y|}{|x-y|}\right)&(1-|y|^p)e^{4u}dy\\
&\ge(-\log|x|+C)\int_{B_{R_0}\setminus B_1}(1-|y|^p)e^{4u}dy\\
&\ge(-\log|x|+C)\int_{B_{R_0}}(1-|y|^p)e^{4u}dy\\
&\ge(-\log|x|+C)(\Lambda+\varepsilon)\ge-(\Lambda+\varepsilon)\log|x|
\end{align*}
where we used the fact that 
$\int_{B_{R_0}\setminus B_1}(1-|y|^p)e^{4u}dy\le\int_{B_{R_0}}(1-|y|^p)e^{4u}dy$.
Concerning the integral over $B_1$ we have
\begin{align*}
\int_{B_1}\log\left(\frac{|x-y|}{|y|}\right)&(1-|y|^p)e^{4u}dy\le\int_{B_1}\log\left(\frac{|x-y|}{|y|}\right)e^{4u}dy\\
&=\int_{B_1}\log\left(\frac{1}{|y|}\right)e^{4u}dy+\int_{B_1}\log(|x-y|)e^{4u}dy\le C
\end{align*}
using Holder's inequality. Therefore, we obtain
\begin{equation}\label{u}
\int_{B_{R_0}}\log\left(\frac{|y|}{|x-y|}\right)(1-|y|^p)e^{4u}dy\ge-(\Lambda+\varepsilon)\log|x|-C.
\end{equation}
We observe that $\log|x-y|\ge0$ for $y\not\in B_1(x)$, $\log|y|\le\log(2|x|)$ for $y\in A_1$, $\int_{A_1}(1-|y|^p)e^{4u} dy\ge-\varepsilon$ and $\log(2|x|)\le2\log|x|$ for $|x|\ge R$, hence we get
\begin{align}\notag
\int_{A_1}\log\left(\frac{|y|}{|x-y|}\right)&(1-|y|^p)e^{4u}dy=\int_{A_1}\log(|y|)\,(1-|y|^p)e^{4u}dy-\int_{A_1}\log(|x-y|)(1-|y|^p)e^{4u}dy\\\notag
&\ge\log(2|x|)\int_{A_1}(1-|y|^p)e^{4u}dy-\int_{B_1(x)}\log(|x-y|)(1-|y|^p) e^{4u}dy\\\notag
&\ge2\log(|x|)\int_{A_1}(1-|y|^p)e^{4u}dy-\int_{B_1(x)}\log(|x-y|)(1-|y|^p) e^{4u}dy\\\label{d}
&\ge-2\varepsilon\log|x|-\int_{B_1(x)}\log(|x-y|)(1-|y|^p) e^{4u}dy.
\end{align}
For $y\in A_2$, in the case $|y|\le2|x|$ we have $\frac{|y|}{|x-y|}\le4$, while in the case $|y|\ge2|x|$ we get $\frac{|y|}{|x-y|}\le2$, so when $y\in A_2$ we have the estimate 
$$\log\left(\frac{|y|}{|x-y|}\right)\le\log4,$$
hence using that $\int_{A_2}(1-|y|^p)e^{4u} dy\ge-\varepsilon$, we obtain
\begin{equation}\label{t}
\int_{A_2}\log\left(\frac{|y|}{|x-y|}\right)(1-|y|^p)e^{4u}dy\ge\log(4)\int_{A_2}(1-|y|^p)e^{4u}dy\ge -\varepsilon\log4.
\end{equation}
Putting together \eqref{u}, \eqref{d} and \eqref{t}, possibly taking $R$ larger, we get
$$v(x)\ge-\frac{1}{8\pi^2}(\Lambda+5\varepsilon)\log|x|+\frac{1}{8\pi^2} \int_{B_1(x)}\log\left(\frac{1}{|x-y|}\right)(1-|y|^p) e^{4u}dy$$
\endproof

\noindent From \eqref{v>} changing signs, it follows that for any $\varepsilon>0$ there is $R>0$ such that for $|x|\ge R$ 
\begin{equation}\label{-vsopra}
-v(x)\le\frac{1}{8\pi^2}(\Lambda+5\varepsilon)\log|x|+\frac{1}{8\pi^2}\int_{B_1(x)}\log(|x-y|)(1-|y|^p) e^{4u}dy.
\end{equation}


\subsection{A Liouville type-theorem}

To prove a Liouville-type theorem (see Theorem \ref{teoliouville} below) we will need the following useful result.

\begin{lemma}\label{lemliouville}
Let $u$ be a measurable function such that $(1-|y|^p)e^{4u}\in L^1(\mathbb{R}^4)$. Then for any $x\in\mathbb{R}^4$ 
$$\fint_{B_r(x)}u^+dy\to 0, \quad\text{as}\,\,r\to+\infty.$$
\end{lemma}

\proof
Let $x\in\mathbb{R}^4$ be fixed, using the fact that $4u^+\le e^{4u}$ we get
\begin{equation}\label{a1}
4\fint_{B_r(x)}u^+dy\le\fint_{B_r(x)}e^{4u}dy=\frac{C}{r^4}\int_{B_r(x)} \frac{1}{1-|y|^p}(1-|y|^p)e^{4u}dy.
\end{equation}
Observing that for $y\in B_r(x)$ we have $|y|\le r+|x|$ and $|y|\ge-r+|x|$, we obtain the following inequalities
$$\frac{1}{1-|y|^p}\le\frac{1}{1-(|x|+r)^p},$$
and 
$$\frac{1}{1-|y|^p}\ge\frac{1}{1-(|x|-r)^p},$$
by means of them we get
$$4\fint_{B_r(x)}u^+dy\le\frac{C}{r^4}\Bigg[\frac{1}{1-(|x|+r)^p}\int\limits_{B_r(x)\cap B_1}(1-|y|^p)e^{4u}dy+\frac{1}{1-(|x|-r)^p}\int\limits_{B_r(x)\cap B_1^c}(1-|y|^p)e^{4u}dy\Bigg]$$
\noindent since $\int_{B_r(x)}(1-|y|^p)\,e^{4u}dy<\infty$, we can estimate \eqref{a1} with $O(r^{-p-4})$ as $r\to\infty$. The claim follows as $r\to+\infty$ since by assumption $p>0$.
\endproof

We are now in the position to prove the following Liouville-type theorem, which will be crucial to prove that $u-v$ is a polynomial of degree at most 2.

\begin{theorem}\label{teoliouville}
Consider $h:\R^4\to\R$ such that
$$\Delta^2 h=0\qquad\text{and}\qquad h\le u-v.$$
Assume that $(1-|y|^p)e^{4u}\in L^1(\R^4)$, $v\in L^1_{loc}(\mathbb{R}^4)$ and further that \eqref{-vsopra} holds. Then $h$ is a polynomial of degree at most $2$.
\end{theorem}

\proof 
We take some ideas from the proof of Theorem 6 in \cite{Mclass}, but this proof is more delicate since our estimate for $-v$ contains the singular integral
\begin{equation}\label{A}
A(x):=\frac{1}{8\pi^2}\int_{B_1(x)}\log(|x-y|)(1-|y|^p) e^{4u}dy. 
\end{equation}
By elliptic estimates for biharmonic functions (see Proposition 4 in \cite{Mclass}) for any $x\in\R^4$ we have
\begin{equation}\label{D3h}
|D^3h(x)|\le\frac{C}{r^3}\fint_{B_r(x)}|h(y)|dy=-\frac{C}{r^3}\fint_{B_r(x)}h(y)dy+\frac{2C}{r^3}\fint_{B_r(x)}h^+dy. 
\end{equation}
From Pizzetti's formula (refer e.g. to\,\,\cite{Mclass}) we have
\begin{equation}\label{pizzetti}
\fint_{B_r(x)}h(y)dy=O(r^2), \quad \text{as}\,\,\,r\to\infty.
\end{equation}
In order to estimate the term $\frac{2C}{r^3}\fint_{B_r(x)}h^+dy$, we observe that
$$\fint_{B_r(x)}h^+dy\le\fint_{B_r(x)}u^+dy+C\fint_{B_r(x)}(-v)^+dy,$$
thanks to Lemma \ref{lemliouville} the term $\fint_{B_r(x)}u^+dy\to0$. 
Using Tonelli's theorem, we can prove that $A\in L^1(\mathbb{R}^4)$ as follows
\begin{align*}
\int_{\R^4} |A(x)|\,dx&=\int_{\R^4}\Big| \frac{1}{8\pi^2}\int_{B_1(x)}\log(|x-y|)(1-|y|^p) e^{4u}dy\Big|dx\\
&\le\int_{\R^4}\frac{1}{8\pi^2}\int_{B_1(x)} \big|\log|x-y|\big|\,\big|1-|y|^p\big| e^{4u} dy\,dx\\
&=C\int_{\R^4}\int_{\R^4} \chi_{|x-y|<1}\big|\log|x-y|\big|\,\, \big|1-|y|^p\big| e^{4u}dy\,dx\\
&= C\int_{\mathbb{R}^4}\big|1-|y|^p\big|e^{4u}\int_{B_1(y)} \log\left(\frac{1}{|x-y|}\right)dx\,dy\\
&=C\int_{\R^4}\big|1-|y|^p\big| e^{4u}dy<\infty.
\end{align*}
Since $\frac{1}{8\pi^2}(\Lambda+5\varepsilon) \log|x|+A(x)\ge0$ for $|x|\ge R>2$, we get
$$(-v)^+\le \frac{1}{8\pi^2}(\Lambda+5\varepsilon)\log|x|+A(x)$$
for $x\in\mathbb{R}^4\setminus B_R$. Taking into account that $A(x)\in L^1$ we obtain
$$\fint_{B_r(x)}(-v)^+dy\le C\fint_{B_r(x)}\log(|y|+1)\,dy+\fint_{ B_r(x)} A(y)dy
\le C\log r+\frac{C}{r^4}$$
(if $y\in B_R(0)$ the previous estimate for $(-v)^+$ does not hold, we can overcome this problem since $v\in L^1_{loc}(\mathbb{R}^4)$). Hence,
\begin{equation}\label{h+}
\fint_{B_r(x)}h^+dy\le\fint_{B_r(x)}u^+dy+\frac{C}{r^4}+C\log r.
\end{equation}
From \eqref{pizzetti} and \eqref{h+}, we get that all terms in \eqref{D3h} go to 0 as $r\to \infty$, hence we obtain $D^3 h\equiv 0$.
\endproof


\section{Proof of Theorem \ref{1.2}}

In order to prove that all solutions to problem \eqref{1}-\eqref{2} have the form $v+P$, where $v$ behaves logarithmically at infinity and $P$ is an upper-bounded polynomial of degree at most 2, we proceed by steps.

\begin{theorem}\label{teopol}
Let $u$ be a solution to problem \eqref{1}-\eqref{2} such that $(1-|y|^p)e^{4u} \in L^1 (\mathbb{R}^4)$ and $v$ as in \eqref{v}. Then 
$$u=v+P$$
where $P$ is a polynomial of degree at most $2$. Moreover, $\Delta u(x)$ can be represented by
\begin{equation}\label{deltau}
\Delta u(x)=-\frac{1}{4\pi^2}\int_{\mathbb{R}^4}\frac{1}{|x-y|^2}(1-|y|^p)e^{4u}dy+C, 
\end{equation}
where $C$ is a constant.
\end{theorem}

\proof
Consider $P=u-v$, we have $\Delta^2 P=0$. From \eqref{-vsopra} using Theorem \ref{teoliouville} we can conclude that $u=v+P$ where $P$ is a polynomial of degree at most $2$. Hence $\Delta u=\Delta v+\Delta P$, it follows immediately that $\Delta P=C$ where $C$ is a constant and from \eqref{v} we have that $\Delta v=-\frac{1}{4\pi^2}\int_{\mathbb{R}^4}\frac{1}{|x-y|^2}(1-|y|^p) e^{4u}dy$.
\endproof

Let us prove that the polynomial $P$ is upper-bounded. 

\begin{prop}\label{sup p}
Let $P$ be the polynomial of Theorem \ref{teopol}. Then 
$$\sup\limits_{x\in\mathbb{R}^4} P(x)<+\infty.$$
\end{prop}

\proof
\textbf{Step 1. Estimate of the term $|A(x)|$.} We take some ideas from the proof of Lemma 13 in \cite{Mclass} and the one of Lemma 2.4 in \cite{Lin}, but our case is more challenging.  In what follows, $C$ denotes a generic constant which may change from line to line and also within the same line.
We observe that
\begin{align*}
|A(x)|&=\bigg|\frac{1}{8\pi^2}\int_{B_1(x)}\log\left(|x-y|\right)(1-|y|^p) e^{4u}dy\bigg|\le\frac{1}{8\pi^2}\int_{B_1(x)}\log\left(\frac{1}{|x-y|}\right)\left|1-|y|^p\right| e^{4u}dy \\
&=\frac{1}{8\pi^2}\int_{B_1(x)\setminus B_\tau(x)}\log\left(\frac{1}{|x-y|}\right)\left |1-|y|^p \right| e^{4u}dy+\frac{1}{8\pi^2}\int_{B_\tau(x)}\log\left(\frac{1}{|x-y|}\right)\left|1-|y|^p\right| e^{4u}dy
\end{align*}
where $\tau\in(0,1)$ will be fixed later. Since $\log\left(\frac{1}{|x-y|}\right)\in(0,-\log\tau)$ for $y\in B_1(x)\setminus B_\tau(x)$ and by assumption $(1-|y|^p)e^{4u}\in L^1(\mathbb{R}^4)$, we have 
$$\int_{B_1(x)\setminus B_\tau(x)}\log\left(\frac{1}{|x-y|}\right)\left |1-|y|^p \right| e^{4u}dy<C.$$ 
In order to estimate the integral over $B_\tau(x)$ we proceed as follows.
By Holder's inequality, we get
\begin{align}\notag
\int_{B_\tau(x)}\log&\left(\frac{1}{|x-y|}\right)\left|1-|y|^p\right| e^{4u}dy \le\left(\int_{B_\tau(x)}\log\left(\frac{1}{|x-y|}\right)^2dy\right)^{1/2}
\left(\int_{B_\tau(x)}\left|1-|y|^p\right|^2 e^{8u}dy\right)^{1/2}\\ \label{holder Btau}
&\le\left(\int_{B_\tau(x)}\log\left(\frac{1}{|x-y|}\right)^2dy\right)^{1/2}
\left(\int_{B_\tau(x)}\left|1-|y|^p\right|^4 dy\right)^{1/4}
\left(\int_{B_\tau(x)} e^{16u}dy\right)^{1/4}.
\end{align}
Fix $0<\varepsilon_0<1$, we can choose $R_0>6$ sufficiently large such that 
\begin{equation}\label{a}
\int_{B_4(x)}\left|1-|y|^p\right|e^{4u}dy\le\varepsilon_0
\end{equation}
for $|x|\ge R_0$. Let $h$ be the solution of 
$$\begin{cases} \Delta^2 h=f\quad\text{on}\,\,B_4(x)\\
h=\Delta h=0\quad \text{on}\,\,\partial B_4(x)\end{cases}$$
where $f(y)=(1-|y|^p)e^{4u(y)}$, then by Theorem 7 in \cite{Mclass} (refer also to Lemma 2.3 \cite{Lin}) for any $k\in\left(0,\frac{8\pi^2}{\|f\|_{L^1(B_4(x))}}\right)$, we have 
$e^{4k|h|}\in L^1(B_4(x))$, and
\begin{equation}\label{e^4ph}
\int_{B_4(x)}e^{4k|h|}dy\le C
\end{equation}
where $C$ is a constant which depends on $k$ but is independent from $x$. For $y\in B_4(x)$ define $q(y):=u(y)-h(y)$, then $q$ satisfies
$$\begin{cases} \Delta^2 q=0\quad\text{on}\,\,B_4(x)\\
\Delta q=\Delta u\quad\text{and}\quad q=u\quad \text{on}\,\,\partial B_4(x)\end{cases}.$$
Integrating equation $\Delta^2 u=(1-|y|^p)e^{4u}$ on $B_\rho(x)$ we get
\begin{equation}\label{23b}
\int_{\partial B_\rho(x)}\frac{\partial}{\partial r}(\Delta u)d\sigma=\int_{B_\rho(x)} (1-|y|^p)e^{4u}dy.    
\end{equation}
Dividing by $\omega_4\rho^3$ and integrating with respect to $\rho$ from 0 to $R$ (we assume $R<5$), using Fubini's theorem, we obtain 
$$\int_0^R\frac{1}{\omega_4 \rho^3}\int_{\partial B_\rho(x)}\frac{\partial }{\partial r}(\Delta u) d\sigma \,d\rho=\fint_{\partial B_R(x)}\Delta u\,d\sigma-\Delta u(x)$$
and similarly
$$\int_0^R\frac{1}{\omega_4\rho^3}\int_{B_\rho(x)}(1-|y|^p)e^{4u}dy\,d\rho=\frac{1}{4\pi^2}\int_{B_R(x)}(1-|y|^p)e^{4u}\left[\frac{1}{|x-y|^2}-\frac{1}{R^2}\right]dy.$$
Hence
$$\fint_{\partial B_R(x)}\Delta u\,d\sigma=\Delta u(x)+\frac{1}{4\pi^2} \int_{B_R(x)}(1-|y|^p)e^{4u}\left[\frac{1}{|x-y|^2}-\frac{1}{R^2}\right]dy$$
by means of identity \eqref{deltau} we get
$$-\fint_{\partial B_R(x)}\Delta u\,d\sigma=\frac{1}{4\pi^2} \int_{|x-y|\ge R}\frac{1-|y|^p}{|x-y|^2}e^{4u}dy+\frac{1}{4\pi^2 R^2}\int_{B_R(x)}(1-|y|^p)e^{4u}dy -C_1.$$
If we take $R=4$ we have
\begin{equation}\label{stima delta u}
-\fint_{\partial B_4(x)}\Delta u\,d\sigma\le C.    
\end{equation}
Let $G$ be the Green's function for the operator $\Delta$ on $B_4(x)$, namely
$$\Delta G=\delta_x \qquad\text{and} \qquad G=0 \,\,\,\text{on}\,\,\,\partial B_4(x),$$
we have 
$$-\Delta q(x)=-\int_{\partial B_4(x)}\frac{\partial G}{\partial n} \Delta u\,d\sigma= -\int_{\partial B_4(x)}c_0\, \Delta u\,d\sigma\le C$$
where by \cite[Lemma 12]{Mclass} $c_0>0$ and in the last inequality we used \eqref{stima delta u}. Since $c_0$ is a positive constant there exist some $\tau\in(0,1)$ such that if $\xi\in B_{4\tau}(x)$ and $G_\xi$ is the Green's function defined by
$$\Delta G_\xi=\delta_\xi, \quad G_\xi=0 \,\,\,\text{on}\,\,\,\partial B_4(x),$$
then $$0\le \frac{\partial G_\xi(\eta)}{\partial r}\le C, \quad \text{for}\,\,\eta\in \partial B_4(x), \quad r:=\frac{\eta-x}{4}$$
and as before we get 
\begin{equation}\label{26}
-\Delta q(y)\le C \quad\text{on}\,\, B_{4\tau}(x).
\end{equation}
Define $\tilde{q}(y):=-\Delta q(y)$, obviously $q$ satisfies 
$$\begin{cases}
\Delta q(y)=-\tilde{q}(y)\quad\text{in}\,\, B_4(x)\\ q=u\quad\text{on}\,\,\partial B_4(x)
\end{cases}$$
hence by elliptic estimates (refer to \cite[Theorem 8.17]{GT}) for any $\ell>1$ and $\sigma>2$ 
$$\sup\limits_{B_\tau(x)}q\le c(\ell,\sigma)\left(\|q^+\|_{L^\ell(B_{4\tau}(x))} +\|\tilde{q}\| _{L^\sigma(B_{4\tau}(x))}\right).$$
From \eqref{26} we get $\|\tilde{q}\|_{L^\sigma(B_{4\tau}(x))}\le C$. 
Since $q=u-h$, it follows that $q^+(y)\le u^+(y)+|h(y)|$ for $y\in B_4(x)$ and hence
$$\int_{B_{4\tau}(x)}(q^+)^2\le C\int_{B_{4\tau}(x)}e^{2q^+}\le C\left(\int_ {B_{4\tau}(x)}e^{4u^+}\right)^{1/2}\left(\int_{B_{4\tau}(x)}e^{4|h|}\right)^{1/2}.$$
Note that
\begin{equation}\label{aa}
e^{4u^+}\le 1+e^{4u}\le1+\left|1-|y|^p\right| e^{4u},\quad \text{for}\,\, |y|\ge2^{1/p}.
\end{equation}
Since $B_{4\tau}(x)\subset B_{2^{1/p}}^c$ (eventually choosing $R_0$ greater) from \eqref{a} we get 
$$\int_{B_{4\tau}(x)}e^{4u^+}dy\le C$$
finally, together with \eqref{e^4ph}, we obtain $$\|q^+\|_{L^2(B_{4\tau}(x))}<C.$$
In this way we have shown that
$$\sup\limits_{B_\tau(x)}q\le C,$$
hence for $y\in B_\tau(x)$ we have
$$u(y)=q(y)+h(y)\le C+|h(y)|,$$
from which we get 
\begin{equation}\label{e^16u}
\int_{B_\tau(x)}e^{16u}dy\le C\int_{B_\tau(x)}e^{16|h|} dy<C.
\end{equation}
where in the last inequality we used \eqref{e^4ph}. From \eqref{holder Btau} using \eqref{e^16u} we get that
$$\int_{B_\tau(x)}\log\left(\frac{1}{|x-y|}\right)\big|1-|y|^p\big| e^{4u}dy\le C$$
hence $|A(x)|\le C$.\\
\textbf{Step 2.} Let us define 
$$f(r):=\sup\limits_{x\in\partial B_r} P(x)$$
and assume by contradiction that $\sup\limits_{\mathbb{R}^4}P=+\infty$.
From Theorem 3.1 in \cite{Go} it must exist $s>0$ such that 
$$\lim\limits_{r\to+\infty}\frac{f(r)}{r^s}=+\infty.$$
Moreover, $P$ is a polynomial of degree at most 2, hence $|\nabla P(x)|\le c |x|$ for $|x|$ large. From Lemma \ref{v>=} using the fact that $|A(x)|\le C$, we get that there is $R$ sufficiently large such that for every $r\ge R$ we can find $x_r$ with $|x_r|=r$ such that 
$$u(y)=v(y)+P(y)\ge r^s, \quad\text{for}\,\,|y-x_r|\le\frac{1}{r}.$$
We consider 
\begin{align*}
-\int_{B_1^c}(1-|y|^p)e^{4u}dy&=\int_{B_1^c}(|y|^p-1)e^{4u}dy\ge\int_{B_R^c}(|y|^p-1)e^{4u}dy\\
&\ge\int_{B_R^c}e^{4u}dy\ge C\int_R^{+\infty}\int_{\partial B_r\cap B_{1/r}(x_r)} e^{4u} d\sigma\,dr\\
&\ge C\int_R^{+\infty}\int_{\partial B_r\cap B_{1/r}(x_r)} e^{4r^s} d\sigma\,dr\ge C\int_R^{+\infty}\frac{e^{4r^s}}{r^3}dr=+\infty.
\end{align*}
which is absurd since $(1-|y|^p)e^{4u}\in L^1(\mathbb{R}^4)$.
\endproof

The fact that $P$ is upper-bounded implies that $P$ has even degree, hence or $P$ has degree $2$ or is constant.

\begin{prop}\label{v>new}
Let $u$ be a solution to \eqref{1}-\eqref{2} such that $(1-|y|^p)e^{4u}\in L^1(\R^4)$ and $v$ as defined in \eqref{v}. Then given any $\varepsilon>0$, there exists $R=R(\varepsilon)$ such that for $|x|\ge R$, $v(x)$ satisfies
\begin{equation}\label{vnew}
v(x)\ge -\frac{1}{8\pi^2}(\Lambda+6\varepsilon)\log|x|.
\end{equation}
Moreover, we have 
\begin{equation}\label{Dvlim}
\lim\limits_{|x|\to+\infty}\Delta v(x)=0.
\end{equation}
\end{prop}

\proof 
First we prove \eqref{vnew}. From Lemma \ref{v>=} for any $\varepsilon>0$ there exists $R>0$ such that for $|x|\ge R$ 
$$v(x)\ge -\frac{1}{8\pi^2}(\Lambda+5\varepsilon)\log|x|-\frac{1}{8\pi^2}\int_{B_1(x)}\log(|x-y|)(1-|y|^p) e^{4u}dy.$$
Notice that a priori the second term on the right-hand side may be very little, but in the proof of Proposition \ref{sup p} we have proven that $|A(x)|\le C$, so \eqref{vnew} follows at once.\\
\indent Now we prove \eqref{Dvlim}. Differentiating we have 
$$\Delta v(x)=-\frac{1}{4\pi^2}\int_{\mathbb{R}^4} \frac{1}{|x-y|^2} (1-|y|^p)e^{4u}dy.$$
For any $\sigma>0$, by dominated convergence 
$$\int\limits_{\mathbb{R}^4\setminus B_\sigma(x)} \frac{(1-|y|^p)e^{4u}}{|x-y|^2}dy\to 0, \quad\text{as}\,\,|x|\to+\infty.$$
By Holder's inequality we get
$$\int\limits_{B_\sigma(x)} \frac{(1-|y|^p)e^{4u}}{|x-y|^2}dy\le 
\left(\int\limits_{B_\sigma(x)} \frac{(1-|y|^p)^k}{|x-y|^{2k}}dy\right)^{1/k}
\left(\int\limits_{B_\sigma(x)}e^{4k'u}dy\right)^{1/k'}$$
if $\sigma$ is small enough, by \eqref{e^16u} we can conclude.
\endproof 

\proof[Proof of Theorem \ref{1.2}] It follows from Lemma \ref{v<}, Theorem \ref{teopol}, Proposition \ref{sup p} and Proposition \ref{v>new}.
\endproof

\begin{corollary}
Any solution $u$ to \eqref{1}-\eqref{2} is bounded from above.
\end{corollary}

\proof
The solution $u$ is continuous and $u=v+P$. Moreover from \eqref{valto} we have that $v(x)\le C$ on $B_4^c$ and from Proposition \ref{sup p} we have $\sup_{\mathbb{R}^4}P(x)<+\infty$.
\endproof


\section{Scalar curvature}

In this section we prove Theorem \ref{1.3}. First, recall that if $u$ is a solution to \eqref{1}-\eqref{2}, then the metric $g_u=e^{2u}g_{\mathbb{R}^4}$ is conformal to the flat metric on $\mathbb{R}^4$ and has $Q$-curvature equal to $1-|x|^p$. If we consider a metric $g_u$ conformal to a flat metric, the conformal change of scalar curvature is given by the following formula
$$R_{g_u}=-6e^{-2u}\left(\Delta u+|\nabla u|^2\right).$$

In the case when $u$ is a normal solution to problem \eqref{1}-\eqref{2}, A. Hyder and L. Martinazzi (refer to Theorem 1.2 in \cite{HM20}) proved that for $\ell=1,2,3$ 
$$|\nabla^\ell u(x)|=O(|x|^{-\ell}),\quad\text{as}\,\,|x|\to+\infty.$$
Moreover, from the fact that $|\Delta u(x)|\to 0$ as $|x|\to+\infty$ and $\Delta^2u\le0$ for $|x|\ge1$, we get $\Delta u(x)\to0$ as $|x|\to+\infty$. Hence, if $u$ normal solution, we have
$$\lim\limits_{|x|\to+\infty}R_{g_u}=0.$$

In the case when $u=v+p$ non-normal solution, we have proved that the polynomial $p$ has degree 2, hence $\Delta p$ is constant and $deg |\nabla p|^2=2$. In this way we obtain that 
$$\limsup\limits_{|x|\to+\infty}(\Delta p+|\nabla p|^2)=+\infty.$$
Differentiating \eqref{v} we get
$$\nabla v(x)=- \frac{1}{8\pi^2}\int_{\mathbb{R}^4}\frac{x-y}{|x-y|^2}(1-|y|^p)e^{4u}dy$$
hence $|\nabla v(x)|\to 0$ as $|x|\to+\infty$; from Proposition \ref{v>new} we get $\lim\limits_{|x|\to+\infty}\Delta v(x)=0$. Then we obtain
$$\limsup\limits_{|x|\to+\infty}\left(\Delta u+|\nabla u|^2\right)= \limsup\limits_{|x|\to+\infty}\left(\Delta p+|\nabla p|^2\right)=+\infty$$
and hence since $e^{-2u}>0$ we get
$$\liminf\limits_{|x|\to+\infty}R_{g_u}=-\infty.$$

\textbf{Acknowledgements}  I would like to thank Professor Luca Martinazzi for suggesting me the problem and for interesting discussions and suggestions.


\vspace{1.5cm}
\flushright{\textsc{Chiara Bernardini}\\
Dipartimento di Matematica, Università di Padova\\
Via Trieste 63, 35121 Padova (Italy)\\
\mail{chiara.bernardini@math.unipd.it}}


\begin{thebibliography}{90}

\bibitem{B} Branson T.; \textit{The functional determinant}, Global Analysis Research Center Lecture Notes Series, no. 4, Seoul National University (1993).
\doi{10.1090/S0002-9939-1991-1050018-8}

\bibitem{B1} Branson T.; \textit{Sharp inequality, the functional determinant and the complementary series}, Trans. Amer. Math. Soc. \textbf{347} (1995), 3671-3742. \doi{10.1090/S0002-9947-1995-1316845-2}

\bibitem{BO} Branson T., {\O}rsted B.; \textit{Explicit functional determinants in four dimensions}, Proc. Amer. Math. Soc.\textbf{113} (1991), 669-682.
\doi{10.2307/2048601}

\bibitem{CC} Chang S.-Y. A., Chen W.; \textit{A note on a class of higher order conformally covariant equations}, Discrete Contin. Dynam. Systems \textbf{7} (2001), no. 2, 275-281. \doi{10.3934/dcds.2001.7.275}


\bibitem{CY} Chang S.-Y. A., Yang P.; \textit{On uniqueness of solutions of $n$-th order differential equations in conformal geometry}, Math. Res. Lett. \textbf{4} (1997), 91-102. \doi{10.4310/MRL.1997.V4.N1.A9}

\bibitem{CL1} Chen W., Li C.; \textit{Classification of solutions of some nonlinear elliptic equations}, Duke Math. J. \textbf{63} (1991), 615-622.
\doi{10.1215/S0012-7094-91-06325-8}

\bibitem{Ch} Chern S.-S.; \textit{A simple intrinsic proof of the Gauss-Bonnet theorem for closed Riemannian manifolds}, Ann. Math. \textbf{45} (1944), 747-752. \doi{10.2307/1969302}

\bibitem{FG} Fefferman C., Graham C.R.; \textit{$Q$-curvature and Poincaré metrics}, Math. Res. Lett. \textbf{9} (2002) 139-151. 
\doi{10.4310/MRL.2002.v9.n2.a2}

\bibitem{FH} Fefferman C., Hirachi K.: \textit{Ambient metric construction of Q-curvature in conformal and CR geometries}, Mathematical Research Letters \textbf{10} (2003), 819-831. \doi{10.4310/MRL.2003.v10.n6.a9}

\bibitem{GT} Gilbarg D., Trudinger N.; \textit{Elliptic partial differential equations of second order}, Reprint of the 1998 edition, Classics in Mathematics, Springer-Verlag, Berlin, 2001. xiv+517 pp. ISBN: 3-540-41160-7.

\bibitem{Go} Gorin E.A.; \textit{Asymptotic properties of polynomials and algebraic functions of several variables}, Russ. Math. Surv. \textbf{16} (1)  (1961), 93-119. \doi{10.1070/RM1961v016n01ABEH004100}

\bibitem{GJMS} Graham C.R., Jenne R., Mason L., Sparling G.; \textit{Conformally invariant powers of the Laplacian, I: existence}, J. London Math. Soc. \textbf{46} no.2 (1992), 557-565. \doi{10.1112/jlms/s2-46.3.557}

\bibitem{H} Hyder A.; \textit{Conformally Euclidean metrics on $\mathbb{R}^n$ with arbitrary total $Q$-curvature}, Analysis \& PDE. \textbf{10} (2017), no. 3, 635-652. \doi{10.2140/apde.2017.10.635}

\bibitem{HM14} Hyder A., Martinazzi L.; \textit{Conformal metrics on $\mathbb{R}^{2m}$ with constant $Q$-curvature, prescribed volume and asymptotic behavior}, Discr. Cont. Dynamical Systems - A 35 (2015), 283-299. \doi{10.3934/dcds.2015.35.283}

\bibitem{HM20} Hyder A., Martinazzi L.; \textit{Normal conformal metrics on $\mathbb{R}^4$ with $Q$-curvature having power-like growth}, Journal of Differential Equations \textbf{301} (2021), 37-72. \doi{10.1016/j.jde.2021.08.014}

\bibitem{Lin} Lin C. S.; \textit{A classification of solutions of conformally invariant fourth order equations in $\mathbb{R}^n$}, Comm. Math. Helv. \textbf{73} (1998), 206-231. \doi{10.1007/s000140050052}

\bibitem{Mclass} Martinazzi L.; \textit{Classification of solutions to the higher order Liouville’s equation on $\mathbb{R}^{2m}$}, Math. Z. \textbf{263} (2009), 307-329. \doi{10.1007/s00209-008-0419-1}

\bibitem{M13} Martinazzi L.; \textit{Conformal metrics on $\mathbb{R}^{2m}$ with constant $Q$-curvature and large volume}, Ann. Inst. H. Poincaré Anal. Non Linéaire \textbf{30}:6 (2013), 969–982. \doi{10.1016/j.anihpc.2012.12.007}

\bibitem{P} Paneitz S.; \textit{A quartic conformally covariant differential operator for arbitrary pseudo-Riemannian manifolds}, SIGMA Symmetry Integrability Geom. Methods Appl. \textbf{4} (2008), Paper 036. Preprint (1983). \doi{10.3842/SIGMA.2008.036}

\bibitem{WY} Wei J., Ye D.; \textit{Nonradial solutions for a conformally invariant fourth order equation in $\mathbb{R}^4$}, Calc. Var. Partial Differential Equations \textbf{32} (2008), no. 3, 373-386.
\doi{10.1007/s00526-007-0145-2}
\end{thebibliography}
\end{document}